\documentclass{ajour}

\begin{document}

\journame{Journal of Number Theory}
\articlenumber{}
\yearofpublication{}
\volume{}
\cccline{}


\authorrunninghead{Auer and Top}
\titlerunninghead{Legendre elliptic curves over finite fields}



\title{Legendre elliptic curves over finite fields}

\author{Roland Auer}

\affil{Vakgroep Wiskunde RuG,
P.O. Box 800,
9700 AV Groningen,
The Netherlands}

\email{auer@math.rug.nl}

\author{Jaap Top}
 
\affil{Vakgroep Wiskunde RuG,
P.O. Box 800,
9700 AV Groningen,
The Netherlands}

\email{top@math.rug.nl}

\abstract{
  We show that every elliptic curve over a finite field of odd
  characteristic whose number of rational points is divisible by 4 is
  isogenous to an elliptic curve in Legendre form, with the sole
  exception of a minimal respectively maximal elliptic curve.
  We also collect some results concerning the supersingular
  Legendre parameters.}


\newcommand {\citep}[2] {\cite[p.~#1]{#2}}
\newcommand {\citepp}[2] {\cite[pp.~#1]{#2}}
\newcommand {\comment}[1] {}
\newcommand {\df}[1] {{\em #1}}
\newcommand {\oneq} {\renewcommand{\theequation}{$*$}}
\newcommand {\qt}[1] {`#1'}
\newcommand {\cf} {cf.\ }

\newcommand {\F} {{\bf F}}
\newcommand {\Fb} {\F_{2^n}}
\newcommand {\Fpp} {\F_{p^2}}
\newcommand {\Fppm}[1] {\F_{p^2}^{*#1}}
\newcommand {\Fqm}[1] {{\F_q^*}^{#1}}
\newcommand {\Fpl}{\F_p\sm\{0,1\}}
\newcommand {\Fql}{\F_q\sm\{0,1\}}
\newcommand {\kl}{k\sm\{0,1\}}
\newcommand {\N} {{\bf N}}
\renewcommand {\P} {{\bf P}}
\newcommand {\Q} {{\bf Q}}
\newcommand {\Z} {{\bf Z}}

\newcommand {\card}[1] {\left|#1\right|}
\newcommand {\frob}{\phi}
\newcommand {\orb}[1] {[#1]}
\newcommand {\la} {\lambda}
\newcommand {\ph} {\varphi}
\newcommand {\set}[2] { \{ #1 : #2 \} }
\newcommand {\sm} {\setminus}
\newcommand {\sseq} {\subseteq}
\newcommand {\ti}[1] {\tilde{#1}} 

\newenvironment{map}[1] 
  {\[#1:\begin{array}{rcl}} 
  {\end{array}\]
  \\[-0.5\baselineskip]  
}

\newenvironment{map*}
  {\[\begin{array}{rcl}} 
  {\end{array}\]
  \\[-0.5\baselineskip]  
}

\renewcommand{\thesubsection}{\arabic{subsection}}
\newcounter{substate}
\renewcommand{\thesubstate}{(\alph{substate})}
\newenvironment{substate}
{
  \begin{list}{\bf\thesubstate}
    {\usecounter{substate} 
     \itemindent0em
     \settowidth\labelwidth{\bf(b)} \labelsep0.5em  
     \leftmargin\labelwidth \addtolength\leftmargin\labelsep
     \topsep0.5ex 
     \itemsep0ex}
}
{\end{list}}

\newenvironment {subproof}
{
  \begin{list}{\bf\thesubstate}
    {\usecounter{substate} 
     \leftmargin0em 
     \settowidth\labelwidth{\bf(b)} \labelsep0.5em  
     \itemindent\labelwidth \addtolength\itemindent\labelsep
     \topsep0.5ex 
     \itemsep0ex}
}
{\end{list}}

\begin{article}

\section{Introduction} 
Throughout this paper, $q>1$ denotes a power of an odd prime number $p$,
and $k$ is a field.
Given two elliptic curves $E/k$ and $E'/k$,
all morphisms from $E$ to $E'$ are understood to be defined over $k$.
In particular, we simply write $\mbox{\rm End}(E)$ for the ring of all 
endomorphisms of $E/k$.
The notation $E\simeq E'$ indicates that $E$ is isomorphic to $E'$, and
$E\sim E'$ means that $E$ and $E'$ are isogenous.
The endomorphism of multiplication by $m\in\Z$ on $E$ is denoted by $[m]$.
In case $k=\F_q$, it is a well known fact (see~\cite{Ta}) that $E\sim E'$ if
and only if $\card{E(\F_q)}=\card{E'(\F_q)}$.
The Frobenius endomorphism on an elliptic curve $E/\F_q$ will be denoted
by $\frob=\frob_q$.

For $\mbox{\rm char}(k)\neq2$ and $\la\in\kl$, the 
\df{Legendre elliptic curve}
$E_\la/k$ is given by the equation $y^2=x(x-1)(x-\la)$.
All its $2$-torsion points are rational.
An arbitrary elliptic curve $E/k$ with this property has an equation
of the form $y^2=x(x-\alpha)(x-\beta)$ with $\alpha,\beta\in k^*$
(after a suitable choice of coordinates).
Investigating the possible transformations (see~\cite[III~\S1]{Si})
yields that $E$ is \df{Legendre isomorphic}
(i.e., isomorphic to a Legendre elliptic curve) if and only if at least
one of $\pm\alpha,\pm\beta,\pm(\alpha-\beta)$ is a square in $k$.
This is always true when $k(\sqrt{-1})$ is algebraically closed or
when $k=\F_q$ with $q\equiv3\bmod4$, but not, e.g., for $k=\F_{13}$,
$\alpha=-2$ and $\beta=5$.
So the next question to ask is whether $E$ is isogenous to
a Legendre elliptic curve, or \df{Legendre isogenous}, for short.
For $k=\F_q$, this can be answered affirmatively, with precisely one
exception, which occurs when $q$ is a square.

\begin{theorem}\label{main}
  Let $E/\F_q$ be an elliptic curve.
  Write $q=r^2$, such that $r\equiv1\bmod4$ when $q$ is a square.
  Then $E$ is Legendre isogenous if and only if
  $\card{E(\F_q)}\in4\Z\sm\{(r+1)^2\}$.
\end{theorem}

A proof of this result will be presented in the next section.
In the third section we collect some results concerning
the \df{supersingular Legendre parameters}, i.e., the values
of $\lambda\in\F_q$ for which $E_{\lambda}$ is supersingular.
Section~4 contains some remarks on the \qt{average} number 
$\card{E_{\lambda}(\F_q)}$ when $\lambda$ ranges over $\Fql$,
and the final section considers an analogue in characteristic 2.

Our motivation for studying this originated from the problem
of finding lower bounds for the maximum number $N_q(3)$ of rational points
on genus $3$ curves over $\F_q$. For example, we have used
the main result of the present paper to prove 
 $N_{3^n}(3)\geq 3^n+6\sqrt{3^n}-11$ for every $n\geq 1$.
We plan to present this and similar results in a forthcoming paper.

\section{Isogenies of Legendre elliptic curves}
In this section we assume the characteristic of $k$ to be different
from~$2$.

\begin{lemma}\label{squarelem}
  Let $E/k$ be given by $y^2=(x-\alpha)(x-\beta)(x-\gamma)$
  with $\alpha,\beta,\gamma\in k$, $\alpha\neq\beta\neq\gamma\neq\alpha$.
  Then $(\gamma,0)\in[2]E(k) \iff \gamma-\alpha,\gamma-\beta\in {k^*}^2$.
\end{lemma}

\begin{proof}
  This is true because the homomorphism (see~\cite[Thm.~1.2]{Sf} and
  \cite[X~\S1]{Si})
  \[
    (x-\alpha,x-\beta,x-\gamma) : 
    E(k) \to k^*/{k^*}^2 \times k^*/{k^*}^2 \times k^*/{k^*}^2
  \]
  has kernel $[2]E(k)$ and sends $(\gamma,0)$ to
  $(\gamma-\alpha,\gamma-\beta,(\gamma-\alpha)(\gamma-\beta))$.
\end{proof}

Given an elliptic curve $E/k$ with Weierstrass equation $y^2=f(x)$
and an element $\alpha\in k^*$, we denote by $E^{(\alpha)}/k$
the elliptic curve with equation $\alpha y^2=f(x)$.
Note that $E\simeq E^{(\alpha)}$ for $\alpha\in {k^*}^2$.
If $E/\F_q$ and $\alpha$ is non-square in $\F_q$, then counting
points by means of the quadratic character on $\F_q$ yields
$\card{E(\F_q)}+\card{E^{(\alpha)}(\F_q)} = 2q+2$.

\begin{lemma}\label{isotwist}
Let $E/k$ be an elliptic curve, and let $\alpha\in k^*$ be non-square.
Suppose $E\simeq E^{(\alpha)}$.
Then $j(E)=1728$ and $k(\sqrt{\alpha})=k(\sqrt{-1})$.
\end{lemma}

\begin{proof}
This can be seen from a calculation using the explicit form of
a possible isomorphism (see~\cite[III~\S1 and Appendix~A]{Si}).
Alternatively, one may use the theory of twisting (\cite[X~\S5]{Si}): 
the condition
$E\simeq E^{(\alpha)}$ implies that the cocycle
$\sigma\mapsto [\sigma(\sqrt{\alpha})/\sqrt{\alpha}]$ is trivial in 
$H^1(\mbox{Gal}(\bar{k}/k),\mbox{\rm Aut}(E\otimes\bar{k}))$,
and hence of the form $\sigma\mapsto\sigma(\varphi)\circ\varphi^{-1}$
for some $\varphi\in \mbox{\rm Aut}(E\otimes\bar{k})$.
Then $\varphi$ has order $4$, from
which the lemma easily follows.
\end{proof}

Let us return to the Legendre elliptic curves $E_\la/k$ with $\la\in \kl$.
It is well known (see~\cite[III~\S1]{Si}) that
$E_\la\otimes\bar{k}\simeq E_\mu\otimes\bar{k}$
if and only if
$\mu\in\orb\la:=\{\la,1-\la,1/\la,1-1/\la,1/(1-\la),\la/(\la-1)\}$,
the orbit of $\la$ under the group generated by the two transformations
$\la\mapsto1/\la$ and $\la\mapsto1-\la$ on $\P^1$.

\begin{proposition}\label{squareprop}
  Let $\la\in\F_q\sm\{0,1,-1,2,1/2\}$. 
  The following conditions are equivalent.
  \begin{substate}
  \item $E_\la\simeq E_\mu$ over $\F_q$ for all $\mu\in\orb\la$.
  \item $-1,\la,1-\la\in\Fqm2$.
  \item $E_\la[4](\F_q)\simeq\Z/4\Z\times\Z/4\Z$.
  \end{substate}
  If $E_\la^{(\alpha)}/\F_q$ is not Legendre isomorphic for some
  $\alpha\in\F_q^*$, then the above conditions are satisfied.
\end{proposition}

\begin{proof}
  Since $-1\notin\orb\la$, we know that $j(E_{\la})\neq 1728$.
  From Lemma~\ref{isotwist} and the isomorphisms
   $E_\la^{(-1)}\simeq E_{1-\la}$,
  $E_\la^{(\la)}\simeq E_{1/\la}$ and
  $E_\la^{(1-\la)}\simeq E_{\la/(\la-1)}$ one concludes
  $\mbox{(a)}\Leftrightarrow\mbox{(b)}$.
  The equivalence of (b) and (c) follows directly from
  Lemma~\ref{squarelem}.
  
  If $E_\la^{(\alpha)}$ is not Legendre isomorphic, then
  $\alpha$ must be non-square, and none of the curves $E_\la^{(-1)}$,
  $E_\la^{(\la)}$ and $E_\la^{(1-\la)}$ is isomorphic to
  $E_\la^{(\alpha)}$ since they are all Legendre isomorphic.
  This implies (again using Lemma~\ref{isotwist}) that $-1,\la,1-\la\in\Fqm2$.
\end{proof}

\begin{proposition}\label{ssgrprop}
  Let $p'=(-1)^{(p-1)/2} p$ and suppose $E_\la/\F_q$ is supersingular.
  Then $\la\in\Fpp$ and $E_\la(\Fpp)\simeq\Z/(p'-1)\Z\times\Z/(p'-1)\Z$.
\end{proposition}

\begin{proof}
  Since $E_\la$ is supersingular, it has $j$-invariant
  $j:=j(E_\la)\in\Fpp$ (\cf~\cite[V~Theorem~3.1]{Si}).
  Hence there exists an elliptic curve $E/\Fpp$ such that
  $E_\la\otimes\bar{\F}_p\simeq E\otimes\bar{\F}_p$.
  Multiplication by $p$ on $E$ is purely inseparable of degree
  $p^2=\deg\frob$ (again~\cite[V~Theorem~3.1]{Si}), and therefore it
  factors as $[p]=\psi\circ\frob$ for some automorphism $\psi$ of $E$.
  Assuming $j\neq0,1728$ for a moment implies $\psi=[\pm1]$, hence
  $\frob=[\pm p]$ and
  $E(\Fpp)=E[1-\frob](\bar{\F}_p)=E[p\pm1](\bar{\F}_p)\simeq(\Z/(p\pm1)\Z)^2$.
  As $p\pm1$ is even, $E$ has an equation $y^2=x(x-\alpha)(x-\beta)$ with
  $\alpha,\beta\in\Fpp^*$.
  From $E^{(\alpha)}\simeq E_{\beta/\alpha}$ we conclude
  $\beta/\alpha\in\orb\la$ and therefore $\la\in\Fpp$.
  For $j\in\{0,1728\}$, the latter fact follows from an easy
  calculation.
  
  Therefore we may consider $E_\la/\Fpp$.
  Comparing the list in~\citep{536}{Wa} (see also \cite{Zhu}) with the
  condition $\card{E_\la(\Fpp)}\in4\Z$ leaves us with the cases
  $\card{E_\la(\Fpp)}=(p\pm1)^2$, so
  $E_\la^{(\alpha)}(\Fpp)\simeq(\Z/(p'-1)\Z)^2$
  for suitable $\alpha\in\Fpp^*$.
  Now $p'-1\equiv 0\bmod 4$, which means in particular that
  $(0,0)\in [2]E_{\la}^{(\alpha)}(\Fpp)$.
  By Lemma~\ref{squarelem}, this implies $-\alpha\in\Fppm2$.
  Hence also $\alpha$ is a square in $\Fpp^*$ and thus 
  $E_\la\simeq E_\la^{(\alpha)}$.
\end{proof}

We are now ready to complete the
\begin{demo}{Proof of Theorem~\ref{main}}
  First of all, $E\sim E_\la$ for some $\la\in\Fql$ implies
  $\card{E(\F_q)}=\card{E_\la(\F_q)}\in4\Z$ since $E_\la(\F_q)$ contains
  the whole $2$-torsion subgroup.
  Moreover, if $q$ is a square and $E_\la$ is supersingular, then
  $\frob=[r]$ on $E_\la/\F_q$ by Proposition~\ref{ssgrprop}, and so
  $\card{E_\la(\F_q)}\neq(r+1)^2$.

  To show the opposite direction, we suppose
  $\card{E(\F_q)}\in4\Z\sm\{(r+1)^2\}$ for the rest of the proof.
  If $E$ does not have all its $2$-torsion rational, then $E(\F_q)$ must
  contain a point $P$ of order $4$.
  Choose $Q\in E[2](\bar{\F}_p)\sm\langle[2]P\rangle$.
  Then $\phi_q(Q)\equiv Q\bmod[2]P$, and so $\ti E = E/\langle[2]P\rangle$ 
  does have rational $2$-torsion $\ti E[2](\F_q) =
   \langle P\bmod[2]P,Q\bmod[2]P \rangle \simeq \Z/2\Z\times\Z/2\Z$,
  generated by the images of $P$ and of $Q$ in $E/\langle[2]P\rangle$.
  We therefore may assume that $E$ is given by an equation
  $y^2=x(x-\alpha)(x-\beta)$ with $\alpha,\beta\in\F_q^*$.
  Hence $E^{(\alpha)}\simeq E_\la$ with $\la:=\beta/\alpha$.

  Let us first assume that $E$ is supersingular.
  Once more investigating the list in~\citep{536}{Wa} yields either
  $\card{E(\F_q)}=q+1$ with non-square $q$, or $\card{E(\F_q)}=(r\pm1)^2$
  with $q$ a square.
  In the first case we have $\card{E(\F_q)}=\card{E^{(\alpha)}(\F_q)}$,
  hence $E\sim E^{(\alpha)}$, which implies the theorem.
  In the second case, we must have $E\simeq E^{(\alpha)}\simeq E_\la$ by
  Proposition~\ref{ssgrprop}.
 
  Now suppose $E$ is ordinary.
  If $E_\la^{(\alpha)}$ is Legendre isogenous, the theorem is proven.
  Otherwise we may assume
  $E_\la[4](\F_q)\simeq\Z/4\Z\times\Z/4\Z$ and $-1\in\Fqm2$
  by Proposition~\ref{squareprop}.
  Using R\"uck's theorem~\cite{Ru}, we conclude that $E_\la\sim E'$
  where $E'/\F_q$ is an elliptic curve with
  $E'[4](\F_q)\simeq\Z/4\Z\times\Z/2\Z$.
  Clearly we can choose the coordinates such that  $E'$ has an equation
  $y^2=x(x-\alpha')(x-\beta')$ with  $\alpha',\beta'\in\F_q^*$
  and $(0,0)\in[2]E'(\F_q)$.
  But then $-\alpha'\in\Fqm2$ by Lemma~\ref{squarelem}.
  Thus $E'\simeq E_{\la'}$ with $\la'=\beta'/\alpha'$,
  and this time $E_{\la'}^{(\alpha)}/\F_q$ is Legendre isomorphic
  according to Proposition~\ref{squareprop}.
\end{demo}

\section{Supersingular Legendre parameters}
{}From Proposition~\ref{ssgrprop}, we see that
not only the supersingular $j$-invariants but even the supersingular Legendre 
parameters are in $\Fpp$.
This is well known; compare \citepp{94, 97}{Dw}.
The proof of Proposition~\ref{ssgrprop} moreover shows that these
supersingular Legendre parameters are squares in $\Fpp^*$.
One can prove an even stronger result which also seems to be
well known. See \cite[Theorem~1.9A]{Brock} for a statement of
the results in this section; Brock's approach is rather different from
the one presented here.

\begin{proposition}\label{sslaprop}
  Let $\la\in\Fql$ such that $E_\la$ is supersingular.
  Then $-\la\in\Fppm8$.
\end{proposition}

\begin{proof}
  By Proposition~\ref{ssgrprop}, we have $\la\in\Fpp$ and 
  $E_\la(\Fpp)\simeq(\Z/(p'-1)\Z)^2$ with
  $p':=(-1)^{(p-1)/2}p \equiv1 \bmod 4$.
  In particular, condition (c) of Proposition~\ref{squareprop} is satisfied
  with $q=p^2$, hence $\lambda\in\Fppm2$.
  Let us fix square roots $\sqrt\la,\,\sqrt{-1}=:i\in\Fpp$.
  By~\cite[III~Example~4.5]{Si},
  $E=E_\la/\langle(0,0)\rangle$ has an equation
  $y^2=x(x+(\sqrt\la+1)^2)(x+(\sqrt\la-1)^2)$,
  so $E_\la\sim E\simeq E_{\hat\la}$ with
  $\hat\la:=\bigl(\frac{\sqrt\la+1}{\sqrt\la-1}\bigr)^2$.
  Because $E_{1-\hat\la}$ is supersingular, too, we can conclude
  $1-\hat\la=\bigl(\frac{2i}{\sqrt\la-1}\bigr)^2\sqrt\la\in\Fppm2$.
  This shows that $\la$ is a
  fourth power in $\Fpp$.
  Applying this result to $1-\hat\la$ instead of $\la$ yields
  \oneq
  \begin{equation}\label{laeq}
    \la,\; \frac{1-\hat\la}{(1+i)^4} =
    \frac{\sqrt\la}{(\sqrt\la-1)^2}\in\Fppm4.
  \end{equation}
  The point
  $P:=(\sqrt\la,i(\la-\sqrt\la))\in E_\la(\Fpp)$ has order $4$, namely
  $[2]P=(0,0)$.
  As in the proof of Lemma~\ref{squarelem}, the group homomorphism
  \[
    (x,x-1,x-\la) :
    E(\Fpp)\to\Fpp^*/\Fppm2\times\Fpp^*/\Fppm2\times\Fpp^*/\Fppm2
  \]
  has kernel $[2]E(\Fpp)$ and sends $P$ to
  $(\sqrt\la,\sqrt\la-1,\sqrt\la-\la)$.
  Together with (\ref{laeq}) we obtain the equivalence
  \begin{eqnarray*}
    -1\in\Fppm8 & \iff & 16 \;|\; p^2-1 = (p'-1)(p'+1) \\
                & \iff & p'\equiv1\bmod8 \iff P\in[2]E_\la(\Fpp) \\ 
                & \iff & \sqrt\la-1\in\Fppm2 \iff \la\in\Fppm8.
  \end{eqnarray*}
  Since we already knew that $\la\in\Fppm4$, the desired result
  drops out.
\end{proof}

Recall from~\cite{De} and~\cite{Ig}
that the supersingular Legendre parameters are exactly the
$m:=(p-1)/2$ distinct roots of the Deuring polynomial
$
  H_p(x)=(-1)^m\sum_{k=0}^m{m\choose k}^2 x^k \in \F_p[x].
$
Thus Proposition~\ref{sslaprop} says that $H_p(-x)$ divides
$x^{(p^2-1)/8}-1$.

Concerning the number
$s_p := \card{\set{\la\in\F_p}{H_p(\la)=0}}$ of supersingular
Legendre parameters in $\F_p$, we have the following.
Write $h(-p)$ for the class number of (the ring of integers in)
$\Q(\sqrt{-p})$.

\begin{proposition}\label{ssFpprop}
The number $s_p$ of supersingular Legendre parameters in $\F_p$ satisfies
  \begin{substate}
  \item $s_p=0$ if and only if $p\equiv1\bmod4$.
  \item $s_3=1$.
  \item If $p\equiv3\bmod 4$ and $p>3$, then $s_p=3h(-p)$.
  \end{substate}
\end{proposition}

In the proof, the following lemma will be used.

\begin{lemma}\label{threelambdas}
Assume $p>3$ and $q\equiv 3\bmod 4$, and let $E/\F_q$ be an elliptic
curve with $j$-invariant $j(E)\neq0$.
If $E$ is Legendre isomorphic, then there are exactly $3$ values of
$\la\in \Fql$ such that $E\simeq E_{\la}$.
\end{lemma}

\begin{proof}
Note that 
$E_{\lambda}^{(-1)}\simeq E_{1-\lambda}$ and
$E_{1/\lambda}^{(-1)}\simeq E_{1-1/\lambda}$ and
$E_{1/(1-\lambda)}^{(-1)}\simeq E_{\lambda/(\lambda-1)}$.
Assume $j:=j(E_{\la})\neq0,1728$ for the moment.
Then $\orb\la$ has $6$ elements and,
using Lemma~\ref{isotwist}, exactly one from each pair
$\{\lambda,1-\la\}$, $\{1/\lambda,1-1/\lambda\}$ and 
$\{1/(1-\lambda),\lambda/(\lambda-1)\}$
yields a curve isomorphic to $E_{\lambda}$ over $\F_q$.

The remaining case $j=1728$ corresponds to $\lambda\in\{-1,2,1/2\}$.
These three values are different since we assume the characteristic to
be $>3$.
The curves $E_{-1}$ and $E_2$ are obviously isomorphic.
Moreover, $E_2^{(2)}\simeq E_{1/2}\simeq E_{1/2}^{(-1)}\simeq E_2^{(-2)}$.
Since $q\equiv 3\bmod 4$, one of $2,-2$ is a square in $\F_q^*$, hence 
$E_{1/2}\simeq E_{2}$, and again we find $3$ values of
$\lambda\in\Fql$ giving the same curve.
\end{proof}

\begin{demo}{Proof of Proposition~\ref{ssFpprop}}
(b) holds because $H_3(x)=-x-1$.
We now assume $p>3$.
Then a supersingular elliptic curve over $\F_p$ has $p+1$ rational points.
For $p\equiv1\bmod4$ this number is not divisible by $4$.
Therefore, $s_p=0$ in this case.
Since $s_p>0$ when $p\equiv 3\bmod 4$
(this follows from $H_p(-1)=0$ for such $p$, or alternatively, from
the fact that $H_p$ has odd degree when $p\equiv 3\bmod 4$, while all
irreducible factors have degree $\leq 2$ by Proposition~\ref{ssgrprop}),
(a) follows.

To prove (c), consider a supersingular elliptic curve $E/\F_p$ with 
$3<p\equiv 3\bmod 4$.
Since $\phi^2=[-p]$ on $E/\F_p$, we have
$\Z[\sqrt{-p}]\simeq\Z[\phi]\sseq\mbox{\rm End}(E)$.
Now $\mbox{\rm End}(E)$ is commutative (since all $\F_p$-endomorphisms by
definition commute with $\phi$ and $\mbox{End}(E\otimes\bar{\F}_p)$ has
rank $4$), and therefore $\mbox{\rm End}(E)\sseq\Z[\frac{1-\sqrt{-p}}{2}]$.
{}From this, one concludes that 
$\mbox{\rm End}(E)\simeq\Z[\frac{1-\sqrt{-p}}{2}]$, the ring of integers
in $\Q(\sqrt{-p})$, precisely when
$1-\phi$ is divisible by $2$ in $\mbox{\rm End}(E)$, which
happens if and only if $\phi$ acts trivially on $E[2](\bar{\F}_p)$,
in other words, if and only if all $2$-torsion on $E$ is
$\F_p$-rational.
In particular, a supersingular $E_{\lambda}/\F_p$ satisfies
$\mbox{\rm End}(E_{\lambda})\simeq\Z[\frac{1-\sqrt{-p}}{2}]$.

Conversely, if an elliptic curve $E/\F_p$ satisfies 
$\mbox{\rm End}(E)\simeq\Z[\frac{1-\sqrt{-p}}{2}]$, then
$E$ is supersingular (the trace of an element of norm $p$ is
divisible by $p$ in the latter ring), and by the argument above,
there are $\alpha,\beta\in\F_p^*$ such that $E$ can be
given by an equation $y^2=x(x-\alpha)(x-\beta)$.
Hence $E\simeq E^{(\alpha)}_{\beta/\alpha}\simeq E^{(-\alpha)}_{1-\beta/\alpha}$.
Since $p\equiv 3\bmod 4$, one of $\alpha,-\alpha$ is a square
in $\F_p^*$, and we conclude that $E$ is Legendre isomorphic.
Moreover, $j(E)\neq 0$ because otherwise 
$p\equiv 2\bmod 3$ by supersingularity, while to have all $2$-torsion
rational one would need $p\equiv 1\bmod 3$.
Hence Lemma~\ref{threelambdas} applies, and we find precisely 
$3$ values of $\lambda\in\Fpl$ for which $E\simeq E_{\lambda}$.

The conclusion is that the number $s_p$ of supersingular values of
$\lambda\in\F_p$ equals $3$ times the number of
$\F_p$-isomorphism classes of elliptic curves  $E/\F_p$ with
$\mbox{\rm End}(E)\simeq\Z[\frac{1-\sqrt{-p}}{2}]$.
By \cite[Thm.~4.5]{Wa} (compare \citep{194}{Sch} where a small
correction is given), the latter number equals $h(-p)$.
\end{demo}

It is known that $h(-p)>\frac{1}{55}\log(p)$ (\citep{232}{GZ},
\citep{321}{Oe}).
Hence, in particular, Proposition~\ref{ssFpprop} implies
that for $p\equiv 3\bmod 4$, the number of $\F_p$-rational zeroes
of $H_p$ tends to infinity when $p\rightarrow\infty$.

\section{Some statistics concerning Legendre elliptic curves}

We will briefly discuss some statistical observations
concerning the numbers $\card{E_\la(\F_q)}$.
First of all, these are integers $\equiv 0\bmod 4$, and by the
Hasse inequality, they lie in the interval $[q+1-2\sqrt{q},q+1+2\sqrt{q}]$.
Moreover, if an integer $N=q+1-t$ in this interval does {\em not} occur
as the number of points of some elliptic curve over $\F_q$, then
$\gcd(t,q)\not =1$ (see \cite[Thm.~4.1]{Wa}; in fact this reference
for given $q$ even precisely describes the remaining at most $5$
values of $t$ with $\gcd(t,q)\not =1$ for which an
elliptic curve over $\F_q$ with $N$ points exists).
It follows that there are roughly $\sqrt{q}(1-1/p)$ numbers
$N\equiv0\bmod 4$ which appear as the number of points of some
elliptic curve over $\F_q$.
By our main theorem, all but at most one of these
appear as the number of points of some $E_\lambda/\F_q$.
Since there are $q-2$ elliptic curves $E_\lambda/\F_q$, this
implies that \qt{on average} there are roughly $p\sqrt{q}/(p-1)$
values of $\lambda\in\Fql$ such that
$\card{E_{\lambda}(\F_q)}$ equals a given occurring $N$.

If the numbers $\card{E_\la(\F_q)}$ had an average of $q+1$ over all
$\la\in\Fql$, then
\[
  S(q) := \sum_{\la\in\Fql} \card{E_\la(\F_q)}
\]
would equal $\bar S(q):=(q-2)(q+1)=q^2-q-2$.
But this is impossible for $q\equiv1\bmod4$ because then
$\bar S(q)\equiv2\bmod4$.

\begin{proposition}
  $S(q)=\bar S(q)+1+(-1)^{(q-1)/2}$.
\end{proposition}

\begin{proof}
This can by shown by naively computing
\begin{eqnarray*}
  \ti S(q)&:=&\card{\set{(x,y,\la)\in\F_q^3}{y^2=x(x-1)(x-\la)}}\\
          & =& 2q + \card{\Fql\times\F_q} = q^2, \\
  S_0(q)&:=&\card{\set{(x,y)\in\F_q^2}{y^2=x^2(x-1)}} \\
        & =& 2+2\card{\Fqm2\sm\{-1\}} = q - (-1)^{(q-1)/2}\;\;
          \mbox{\rm  and} \\
  S_1(q)&:=&\card{\set{(x,y)\in\F_q^2}{y^2=x(x-1)^2}}
          = 2+2\card{\Fqm2\sm\{1\}} = q - 1.
\end{eqnarray*}
Then $S(q) = q-2+\ti S(q)-S_0(q)-S_1(q) = q^2-q-1+(-1)^{(q-1)/2}$.
\end{proof}

An alternative method for computing $S(q)$ is by
considering the rational elliptic surface $X\rightarrow \P^1$
corresponding to the Legendre family over the $\lambda$-line.
Compare \citep{56}{GT} for similar calculations.
The surface $X$ has fibre $X_{\la}=E_{\lambda}$ over $\la\in\Fql$.
Over $\lambda=1$ the fibre $X_1$ consists of two $\P^1$'s meeting in two
rational points.
Hence $\card{X_1(\F_q)}=2q$.
Over $\lambda=0$ the fibre $X_0$ also consists of two copies of $\P^1$
meeting in two points; however, these points are rational precisely when
$-1$ is a square in $\F_q$.
This implies $\card{X_0(\F_q)}=2q+1-(-1)^{(q-1)/2}$.
Finally, the fibre $X_{\infty}$ is of Kodaira type $I_2^*$ and
$\card{X_{\infty}(\F_q)}=7q+1$.
The Lefschetz trace formula now shows that $\card{X(\F_q)}=q^2+10q+1$
and hence $S(q) = 
 \card{X(\F_q)}-\card{X_0(\F_q)}-\card{X_1(\F_q)}-\card{X_{\infty}(\F_q)}
= q^2-q-1+(-1)^{(q-1)/2}.$

\section{An analogue in characteristic two}
For the sake of completeness, we consider the situation in
characteristic~$2$.
Let $n\in\N$.
For each $\la\in\Fb^*$, we have the elliptic curve $E_\la/\Fb$ given by
the equation $y^2+xy=x^3+\la$.
Since $j(E_\la)=1/\la$, they are mutually non-isomorphic.

\begin{proposition}
  An elliptic curve $E/\Fb$ satisfies $\card{E(\Fb)}\in4\Z$ if and only
  if $E\simeq E_\la$ for some $\la\in\Fb^*$.
\end{proposition}

\begin{proof}
Recall that $E/\Fb$ is ordinary, i.e., $\card{E(\Fb)}\in2\Z$, 
if and only if
(after a suitable choice of coordinates) it has an equation
$y^2+xy=x^3+\beta x^2+\la$ with $\beta\in\Fb$ and $\la\in\Fb^*$,
and then $j(E)=1/\la$ (see~\cite[Appendix A]{Si}).
Thus we may assume that $E$ has such an equation.
For $\alpha\in\Fb$, we denote by $E^{(\alpha)}$ the elliptic curve with
equation $y^2+xy=x^3+(\alpha+\beta)x^2+\la$.
Then $E\simeq E^{(\alpha)}$ if and only if $\mbox{\rm Tr}(\alpha)=0$, where 
$\mbox{\rm Tr}$
denotes the trace from $\Fb$ to $\F_2$.
Otherwise $E^{(\alpha)}$ is a quadratic twist of $E$ and
$\card{E(\Fb)}+\card{E^{(\alpha)}(\Fb)}=2^{n+1}+2\equiv 2\bmod4$.
It therefore remains to verify that $E_\la(\Fb)\in4\Z$.
Treating the point at infinity and $(0,\sqrt\la)\in E_\la(\Fb)$ separately,
and dividing the equation by $x^2$, we obtain $\card{E_\la(\Fb)}=2+2N$
with
\begin{eqnarray*}
  N &=& \card{\set{x\in\Fb^*}{\mbox{\rm Tr}(x+\la/x^2)=0}} \\
    &=& \card{\set{x\in\Fb^*}{\mbox{\rm Tr}(x)=\mbox{\rm Tr}(\sqrt\la/x)}},
\end{eqnarray*}
which is odd because $x\mapsto\sqrt\la/x$ is an involution on $\Fb^*$
with precisely one fixed point.
\end{proof}

Applying the Frobenius isogeny $\phi_2$ to $E_\lambda$ results 
in the curve $E_{\la^2}$.
Putting $\xi=x$ and $\eta=y+\lambda$, one finds that $E_{\la^2}$ can be
given by the equation $\eta^2+\xi\eta=\xi^3+\la\xi$.
For this equation, a result like the one given above can be found in a
paper by Schoof and van der Vlugt~\citep{172}{SV}.

\begin{acknowledgments}
It is our pleasure to thank Robert Carls, Marius van der Put,
Jasper Scholten and Bart de Smit 
for their interest in this work, and Brad Brock for pointing out some
relevant references.
\end{acknowledgments}

\end{article}
\end{document}